\documentclass{conm-p-l} 
\usepackage{amsrefs,amssymb,enumitem,mathrsfs,xargs}

\usepackage[utf8]{inputenc} 
\usepackage[T1]{fontenc}
\usepackage[english]{babel}
\usepackage[hidelinks]{hyperref}

\usepackage{tikz}
\usetikzlibrary{matrix,arrows}
\usepackage{tikz-cd}

\usepackage[overload]{textcase}
\usepackage{mathtools}

\theoremstyle{plain}

\theoremstyle{definition}

\theoremstyle{remark}

\newif\ifTNS 
\TNStrue

\def\printtheoremname#1{\csname#1name\endcsname}
\def\printtheoremnames#1{\csname#1names\endcsname}

\def\thmref#1#2{\printtheoremname{#1}\ifTNS~\fi\ref{#1:#2}}

\def\uc#1#2{\MakeUppercase{#1}{#2}} 

\newcommand{\DefTheorem}[2]{\newenvironmentx{#1}[2][1=\empty,2=\empty]{%
    \ignorespaces%
    \ifx##2\empty%
      \begin{#2}%
    \else%
      \begin{#2}[{\uc##2}]%
    \fi%
    \ifx##1\empty%
      {}%
    \else%
      \label{#1:##1}%
    \fi%
    \ignorespaces}{\end{#2}\ignorespacesafterend}}

\newcommand{\prfof}[2]{\protect{Proof of~\thmref{#1}{#2}}}

\makeatother

\DefTheorem{Th}{theorem}
\DefTheorem{Prop}{proposition}
\DefTheorem{Cor}{corollary}
\DefTheorem{Lem}{lemma}
\DefTheorem{Def}{definition}
\DefTheorem{Rem}{remark}

\numberwithin{equation}{section}

\newcommand\Define[2][\empty]{\ignorespaces%
  \emph{#2}}%

\tikzset{
  commutative diagrams/.cd,
  arrow style=tikz,
  diagrams={>=stealth},
}

\newcommand\ger{\mathfrak}
\newcommand\sh{\mathcal}

\newcommand\DMO{\DeclareMathOperator}

\newcommand\ev{{\bar 0}}
\newcommand\odd{{\bar 1}}

\newcommand{\defi}{\coloneqq}     

\newcommand{\Fa}{For all }
\newcommand{\fa}{for all }
\newcommand{\scth}{such that }

\newcommand\ie{\emph{i.e.}~}
\newcommand\via{\emph{via}~}

\newcommand\vphi{\varphi}
\newcommand\vrho{\varrho}
\newcommand\eps{\varepsilon}
\newcommand\nats{\mathbb{N}}
\newcommand\ints{\mathbb{Z}}
\newcommand\cplxs{\mathbb{C}}
\newcommand\aff{{\mathbb A}}
\newcommand\sle{\leqslant}
\newcommand\sge{\geqslant}

\DMO\Ad{\mathrm{Ad}}
\DMO\ad{\mathrm{ad}}
\DMO\GL{\mathrm{GL}}
\DMO\id{\mathrm{id}}
\DMO\tr{\mathrm{tr}}
\DMO\diag{\mathrm{diag}}

\makeatletter
\newcommand\Size[7][1]{
  \ifx#20%
    \def\r@l{}\def\r@m{}\def\r@r{}%
  \else%
    \ifx#21%
     \def\r@l{\bigl}\def\r@r{\bigr}\def\r@m{\bigm}%
    \else%
     \ifx#22%
       \def\r@l{\Bigl}\def\r@r{\Bigr}\def\r@m{\Bigm}%
      \else%
       \ifx#23%
          \def\r@l{\biggl}\def\r@r{\biggr}\def\r@m{\biggm}%
        \else
          \ifx#24%
            \def\r@l{\Biggl}\def\r@r{\Biggr}\def\r@m{\Biggm}%
          \fi%
        \fi%
      \fi%
    \fi%
  \fi%
  \ifx#10%
   \def\r@m{}%
  \fi%
  \r@l#3{#4}\r@m#5{#6}\r@r#7%
}%

\makeatother

\makeatletter
\def\Set@Scallop[#1]#2#3{{#1}\Parens{#2}{#3}}
\newcommand\DeclareScalableOperator[2]{%
  \expandafter\def\csname#1\endcsname{\@ifnextchar[{{#2}\Set@Scallop}{{#2}\Set@Scallop[{}]}}
}
\makeatother

\def\DSO{\DeclareScalableOperator}

\DSO{Hom}{\mathrm{Hom}} 
\DSO{End}{\mathrm{End}} 
\DSO{Aut}{\mathrm{Aut}} 
\DSO{GHom}{\underline{\mathrm{Hom}}} 
\DSO{GEnd}{\underline{\mathrm{End}}} 

\newcommand\Set[3]{
  \Size{#1}{\{}{#2}{|}{#3}{\}}%
}%
\newcommand\Dual[3]{
  \Size[0]{#1}{\langle}{#2}{,}{#3}{\rangle}%
}%
\newcommand\Parens[2]{
  \Size[0]{#1}{(}{#2}{}{}{)}
}
\newcommand\Bracks[2]{
  \Size[0]{#1}{[}{#2}{}{}{]}
}
\newcommand\Abs[2]{
  \Size[0]{#1}{\lvert}{#2}{}{}{\rvert}
}
\makeatletter

\makeatletter

\newif\if@smallmat
\newif\if@none
\newif\if@paren
\newif\if@brack
\newif\if@brace
\newif\if@vline
\newenvironment{Matrix}[2][1]
{\ifx#20%
  \@smallmattrue%
\else%
  \@smallmatfalse
\fi%
\ifx#11%
   \@nonefalse\@parentrue\@brackfalse\@bracefalse\@vlinefalse%
\else%
  \ifx#12%
    \@nonefalse\@parenfalse\@bracktrue\@bracefalse\@vlinefalse%
  \else%
    \ifx#13%
     \@nonefalse\@parenfalse\@brackfalse\@bracetrue\@vlinefalse%
    \else%
      \ifx#14%
        \@nonefalse\@parenfalse\@brackfalse\@bracefalse\@vlinetrue%
      \else%
        \ifx#15%
          \@nonefalse\@parenfalse\@brackfalse\@bracefalse\@vlinefalse%
        \else%
          \@nonetrue\@parenfalse\@brackfalse\@bracefalse\@vlinefalse%
        \fi%
      \fi%
    \fi%
  \fi%
 \fi%
 \if@smallmat%
  \if@none%
   \begin{smallmatrix}%
  \else%
    \if@paren%
      \left(\begin{smallmatrix}%
    \else%
      \if@brack%
        \left[\begin{smallmatrix}%
      \else%
        \if@brace%
          \left\{\begin{smallmatrix}%
        \else%
          \if@vline%
            \left\lvert\begin{smallmatrix}%
          \else%
            \left\lVert\begin{smallmatrix}%
          \fi%
        \fi%
      \fi%
    \fi%
  \fi%
 \else%
  \if@none%
    \begin{matrix}%
  \else%
    \if@paren%
      \begin{pmatrix}%
    \else%
      \if@brack%
        \begin{bmatrix}%
      \else%
        \if@brace%
          \begin{Bmatrix}%
        \else%
          \if@vline%
            \begin{vmatrix}%
          \else%
            \begin{Vmatrix}%
          \fi%
        \fi%
      \fi%
    \fi%
  \fi%
 \fi}%
{\if@smallmat%
  \if@none%
   \end{smallmatrix}%
  \else%
    \if@paren%
      \end{smallmatrix}\right)%
    \else%
      \if@brack%
        \end{smallmatrix}\right]%
      \else%
        \if@brace%
         \end{smallmatrix}\right\}%
        \else%
          \if@vline%
            \end{smallmatrix}\right\rvert%
          \else%
            \end{smallmatrix}\right\rVert%
          \fi%
        \fi%
      \fi%
    \fi%
  \fi%
 \else%
  \if@none%
    \end{matrix}%
  \else%
    \if@paren%
      \end{pmatrix}%
    \else%
      \if@brack%
        \end{bmatrix}%
      \else%
        \if@brace%
          \end{Bmatrix}%
        \else%
          \if@vline%
            \end{vmatrix}%
          \else%
            \end{Vmatrix}%
          \fi%
        \fi%
      \fi%
    \fi%
  \fi%
 \fi}%
                         
\makeatother

\begin{document}
  \title[Schur $Q$-functions and Capelli eigenvalues for \lowercase{$\ger q(n)$}]{Schur $Q$-functions and the Capelli eigenvalue problem for the Lie superalgebra \lowercase{$\ger q(n)$}}

  \author[Alldridge]
  {Alexander Alldridge}

  \address{%
    Mathematisches Institut\\
    Mathematisch-Naturwissenschaftliche Fakult\"at\\
    Universit\"at zu K\"oln\\
    Weyertal 86--90\\
    50931 K\"oln\\
    Germany%
  }
  
  \curraddr{%
    Department of Mathematics\\
    University of California at Berkeley\\
    969 Evans Hall\\
    Berkeley, CA 94720\\
    USA\\ 
    \href{https://aalldridge.github.io}{\url{https://aalldridge.github.io}}
  }  

  \email{alldridg@math.uni-koeln.de, aalldridge@berkeley.edu}

  \author[Sahi]
  {Siddhartha Sahi}

  \address{%
    Department of Mathematics\\    
    Rutgers University\\
    110 Frelinghuysen Rd\\ 
    Piscataway, NJ 08854-8019\\
    United States of America
  }
  \email{sahi@math.rutgers.edu}

  \author[Salmasian]
  {Hadi Salmasian}

  \address{%
    Department of Mathematics and Statistics\\
    University of Ottawa\\
    585 King Edward Ave\\
    Ottawa, Ontario\\
    Canada K1N 6N5%
  }

  \email{hadi.salmasian@uottawa.ca}

  \thanks{Alexander Alldridge gratefully acknowledges support by the German Research Council (Deutsche Forschungsgemeinschaft DFG), grant nos.~AL 698/3-1 and ZI 513/2-1, and the Institutional Strategy of the University of Cologne in the Excellence Initiative. The research of Siddhartha Sahi was partially supported by a Simons Foundation grant (509766) and of Hadi Salmasian by an NSERC Discovery Grant (RGPIN-2013-355464). This work was initiated during the Workshop on Hecke Algebras and Lie Theory, which was held at the University of Ottawa. Hadi Salmasian and Siddhartha Sahi thank the National Science Foundation (DMS-162350), the Fields Institute, and the University of Ottawa for funding this workshop}



  \begin{abstract}
    Let $\ger l\defi \ger q(n)\times\ger q(n)$, where $\ger q(n)$ denotes the queer Lie superalgebra. The associative superalgebra $V$ of type $Q(n)$ has a left and right action of $\ger q(n)$, and hence is equipped with a canonical $\ger l$-module structure. We consider a distinguished basis $\{D_\lambda\}$ of the algebra of $\ger l$-invariant super-polynomial differential operators on $V$, which is indexed by strict partitions of length at most $n$. We show that the spectrum of the operator $D_\lambda$, when it acts on the algebra $\mathscr P(V)$ of super-polynomials on $V$, is given by the factorial Schur $Q$-functions of Okounkov and Ivanov. As an application, we show that the radial projections  of the spherical super-polynomials (corresponding to the diagonal symmetric pair $(\ger l,\ger m)$, where $\ger m\defi\ger q(n)$) of irreducible $\ger l$-submodules of $\mathscr P(V)$ are the classical Schur $Q$-functions. As a further application, we compute the Harish-Chandra images of the Nazarov basis $\{C_\lambda\}$ of the centre of  $\ger U(\ger q(n))$.
  \end{abstract}

  \maketitle

  \section{Introduction}\label{SEcIntro}

  Let $G/K$ be a Hermitian symmetric space of tube type. The \emph{Shilov boundary} of $G/K$ is of the form $G/P=K/M$, where $P=LN$ is the Siegel parabolic subgroup and $M=L\cap K$ is a symmetric subgroup of both $K$ and $L$. Let $\ger l$, $\ger m$, and $\ger n$ be the complexified Lie algebras of $L$, $M$, and $N$, respectively. We set $V\defi\ger n$ and regard $V$ as an $L$-module.  In this setting, $V$ has the structure of a simple Jordan algebra.
 
  The polynomial algebra $\mathscr P(V)$  decomposes as the multiplicity-free direct sum of simple $L$-modules $V_\lambda$, indexed naturally by partitions $\lambda$. In this situation one has canonical invariant ``Capelli'' differential operators of the form $\varphi^k \partial(\varphi)^k$, where $\varphi$ is the Jordan norm polynomial. The spectrum of these operators was computed by Kostant and Sahi \cites{ks91,ks93}, and a close connection with reducibility and composition factors of degenerate principal series was established by Sahi \cites{SahiCompositio,SahiCrelle,SahiSh}.

  Sahi showed \cite{s94} that the decomposition of $\mathscr P(V)$ in fact yields a distinguished basis $\{ D_\lambda\}$, called the \Define{Capelli basis}, of the subalgebra of $L$-invariant elements of the algebra $\mathscr{PD}(V)$ of differential operators on $V$ with polynomial coefficients. Moreover, there is a polynomial $c_\lambda$, uniquely characterized by its degree, symmetry, and vanishing properties, such that $D_\lambda$ acts  on each simple summand $V_\mu$ by the scalar $c_\lambda(\mu)$. The problem of characterizing the spectrum of the operators $D_\lambda$ is  referred to as the \emph{Capelli eigenvalue problem}.

  In fact, Sahi \cite{s94} introduced a universal multi-parameter family of inhomogeneous polynomials that serve as a common generalization of the spectral polynomials $c_\lambda$ across all Hermitian symmetric spaces of rank $n$. Later, Knop and Sahi \cite{ks96} studied a one-parameter subfamily of these polynomials, which already contains all the spectral polynomials. They showed that these polynomials  are eigenfunctions of a class of difference operators extending the Debiard--Sekiguchi differential operators. It follows that the top degree terms of the Knop--Sahi polynomials are Jack polynomials, which for special choices of the parameter become spherical functions. 

  These polynomials were later studied from a different point of view by Okounkov and Olshanski, who referred to them as shifted Jack polynomials.

  Subsequently, supersymmetric analogs of the Knop--Sahi shifted Jack polynomials were constructed by Sergeev and Veselov in \cite{SerVes}. More recently, two of us (Sahi and Salmasian \cite{sahi-salmasian}) have extended this circle of ideas  to the setting of the triples $(\ger l, \ger m,V)$ of the form
  \begin{equation}\label{glglosp}
    \begin{gathered}
      (\ger{gl}(m|n)\times \ger{gl}(m|n),\ger{gl}(m|n),\mathrm{Mat}_{m|n}(\cplxs)),\\
      (\ger{gl}(m|2n),\ger{osp}(m|2n),S^2(\cplxs^{m|2n})).
    \end{gathered}
  \end{equation}
  In each of these situations one has, once again, a canonical Capelli basis of differential operators, and \cite{sahi-salmasian} establishes a precise connection to the abstract Capelli problem of Howe and Umeda \cite{HoweUmeda}. It is further shown in Ref.~\cite{sahi-salmasian} that the spectrum of the Capelli basis is given by specialisations of super analogues of Knop--Sahi polynomials, defined earlier by Sergeev and Veselov \cite{SerVes}. In the case of the triple $(\ger{gl}(m|n)\times \ger{gl}(m|n),\ger{gl}(m|n),\mathrm{Mat}_{m|n}(\cplxs))$, these results follow from earlier work of Molev \cite{MolevFSY}, however the case  $(\ger{gl}(m|2n),\ger{osp}(m|2n),S^2(\cplxs^{m|2n}))$ is harder and requires new ideas.

  The Lie superalgebras $\mathfrak{gl}(m|n)$ and $\mathfrak{osp}(m|2n)$ are examples of basic classical Lie superalgebras. Such an algebra admits an even non-degenerate invariant bilinear form and an even Cartan subalgebra, and many results for ordinary Lie algebras extend to this setting, see for instance Ref.~\cite{as}, where spherical representations for the corresponding symmetric pairs are studied. In this paper, we show that the ideas of Ref.~\cite{sahi-salmasian} can actually be extended to \emph{non-basic} Lie superalgebras. More precisely, we consider the case of the \emph{queer Lie superalgebra }$\mathfrak{q}(n)$, usually defined as the subalgebra of $\mathfrak{gl}(n|n)$ of matrices commuting with an odd involution \cite{KacLSA}. For the present purposes, it is convenient to work with a slightly different realization of $\mathfrak{q}(n)$, which we describe below.

  Let $\mathcal{E}$ be the $\mathbb{C}$-algebra generated by an odd element $\varepsilon $, with $\varepsilon ^{2}=1$; thus as a superspace, $\mathcal{E}\cong \mathbb{C}^{1|1}\cong \mathbb{C}\oplus \mathbb{C}\varepsilon $. Let $\mathscr A$ be the  associative superalgebra of $n\times n$ matrices with entries in $\mathcal{E}$. Then $\mathscr A$ is the associative superalgebra of type $Q(n)$, and $\mathfrak{q}(n)$ is isomorphic to $\mathscr A$ regarded as a Lie superalgebra via the graded commutator%
  \[
    [x,y]\defi xy-(-1)^{\Abs0x\Abs0y}yx.
  \]
  In fact $\mathscr A$ is also a Jordan superalgebra via the graded anticommutator, and an $\mathscr A$-bimodule via left and right multiplication. This bimodule structure induces a $\mathfrak{q}(n)\times \mathfrak{q}(n)$-module structure on $V\defi\mathscr A$.

  In this paper, we consider the Capelli eigenvalue problem for the ``diagonal'' triple 
  \begin{equation}\label{qqq}
    (\ger l,\ger m,V)\defi\Parens1{\mathfrak{q}(n)\times \mathfrak{q}(n),\mathfrak{q}(n),\mathscr A}.
  \end{equation}
  We establish a close connection with the Schur $Q$-functions $Q_{\lambda }$ and their inhomogeneous analogues, the \emph{factorial Schur $Q$-functions} $Q_{\lambda }^{\ast }$, which were originally defined by Okounkov and studied by Ivanov \cite{ivanov97}. Our main results are as follows. From Ref.~\cite{cw}, it is known that the space $\mathscr P(V)$ of super-polynomials on $V$ decomposes as a multiplicity-free direct sum of certain $\ger l$-modules $V_{\lambda }$, which are parametrised by strict partitions $\lambda $ of length at most $n$. It follows that $\mathscr P(V^*)$ decomposes as a direct sum of the contragredient $\ger l$-modules $V_\lambda^*$. In Section \ref{Sectstar}, we describe a certain even linear slice $\ger t^*$ to the $M$-orbits on $V^*$. If $p$ is an $\ger {m}$-invariant super-polynomial on $V^*$, then it is uniquely determined by its restriction to $\ger t^*$. This restriction is an ordinary polynomial, and we call it the \emph{$\ger{m}$-radial part of $p.$}

  \begin{Th}[{th1.1}]
    For every $\lambda$, the $\ger l$-module $V_{\lambda}^*$ contains an $\ger{m}$-spherical super-poly\-no\-mial $p_{\lambda }^*$, which is unique up to a scalar multiple. Moreover, up to a scalar, the $\ger m$-radial part of $p_{\lambda }^*$ is the Schur $Q$-function $Q_{\lambda }$.
  \end{Th}
  This is proved in \thmref{Th}{sphpol-q-hompart} below. Now consider the algebra $\mathscr{PD}(V) $ of polynomial coefficient differential operators on $V$. It has an $\ger{l}$-module decomposition
  \[
    \mathscr{PD}(V) \cong 
    \bigoplus_{\lambda ,\mu }V_{\mu }\otimes V_{\lambda
    }^{\ast }\cong \bigoplus _{\lambda ,\mu }\GHom[_{\ger l}]0{V_{\lambda },V_{\mu
        }},
  \]%
  and we write $D_\lambda$ for the differential operator corresponding to the identity map ${\id}_{V_\lambda}\in\GHom[_{\ger l}]0{V_{\lambda },V_{\mu}}$. The $D_{\lambda }$ are the Capelli operators, and they form a basis for the $\ger{l}$-invariant differential operators acting on $\mathscr{P}(V)$. The operator $D_{\lambda }$ acts on each irreducible component $V_{\mu }$ of $\mathscr P(V)$ by a scalar eigenvalue $c_{\lambda}(\mu)$.

  \begin{Th}[{th1.2}]
    The eigenvalues of the Capelli operator $D_\lambda$ are given by the factorial Schur $Q$-function $Q_\lambda^{\ast }$. More precisely, for all $\lambda,\mu$, we have%
    \[
      c_{\lambda }(\mu )=\frac{Q_{\lambda }^{\ast }(\mu) }{Q_{\lambda }^{\ast }(\lambda) }.
    \]
  \end{Th}
  In fact, we prove \thmref{Th}{th1.2} \emph{first} (see  \thmref{Th}{evpol-q} below) and then use it to prove \thmref{Th}{th1.1}.

  Compared to the cases considered in Equation \eqref{glglosp}, the situation in Equation \eqref{qqq} is more complicated. First, since the Cartan subalgebra of $\mathfrak{q}(n)$ is not purely even, the highest weight space of an irreducible finite dimensional $\ger q(n)$-module is not necessarily one-dimensional. Second, unlike the basic classical cases, the tensor product of two irreducible $\mathfrak{q}(n)$-modules is not necessarily an irreducible $\ger q(n)\times \ger q(n)$-module, and sometimes decomposes as a direct sum of two modules which are isomorphic up to parity change. Third, the $\ger m$-spherical vectors in $\mathscr{P}(V)$ are purely odd, whereas the $\ger m$-spherical vectors in $\mathscr{P}(V^*)$ are purely even. These issues add to the difficulties that arise in the proofs in the case of the symmetric pair in Equation \eqref{qqq}.

  In \cite{sergeev}*{Theorem 3}, Sergeev introduced the $\ger q(n)$-analogue of the Harish-Chandra isomorphism
  \begin{equation}\label{HCeta}
    \eta:\mathcal Z(\ger q(n))\longrightarrow \mathscr{P}(\ger h_\ev),
  \end{equation}
  where $\ger h_\ev$ denotes the even part of the Cartan subalgebra of $\ger q(n)$. 
  \thmref{Th}{th1.2} can be reformulated in terms of the map $\eta$, as follows. 
  The image of $\eta$ can be naturally identified with the space of $n$-variable $Q$-symmetric polynomials 
  (see Section \ref{Sec3.2}). 
  We denote the actions of the first and second factors of $\ger l=\ger q(n)\times\ger q(n)$ on $\mathscr P(V)$
  by $L$ and $R$, respectively. Since the typical ``{$\vrho$}-shift'' for the Sergeev--Harish-Chandra isomorphism is equal to zero, we obtain the following reformulation of \thmref{Th}{th1.2}.

  \begin{Th}[{corrz}]
    For every Capelli operator $D_\lambda$, there exists a unique central element $z_\lambda\in\mathcal Z(\ger q(n))$ such that $L(z_\lambda)=D_\lambda$. Furthermore,
    \[
      \eta(z_\lambda)
      (\mu)=\frac{Q_\lambda^*(\mu)}{Q_\lambda^*(\lambda)}.
    \]
  \end{Th}

  \noindent
  The setting of the present paper was also considered by Nazarov, who constructed \cite{nazarov}*{Eq. (4.7)} a family of invariant differential operators $\{I_\lambda\}$ using characters of the Sergeev algebra \cite{sergeev}. Nazarov also defined \cite{nazarov}*{Eq. (4.6)} certain explicit ``Capelli'' elements $\{C_\lambda\}$ in $\sh Z(\mathfrak{q}(n))$, and proved \cite{nazarov}*{Cor. 4.6} that $I_\lambda=\gamma(C_\lambda)$, where $\gamma$ is the left action of $\mathfrak{q}(n)$ on $V$.

  Although our operators $\{D_\lambda\}$ and central elements $\{z_\lambda\}$ are \emph{different} from the $\{I_\lambda\}$ and $\{C_\lambda\}$ defined by Nazarov, one can make an \emph{a posteriori} connection using our \thmref{Prop}{PrptiQlQlstar} below. This allows us to  compute the Harish-Chandra image of Nazarov's central elements $\{C_\lambda\}$. The following result follows immediately from \thmref{Th}{cor1.3}.
  \begin{Th}[thcr]
    The Harish-Chandra image of the operator $C_\lambda$ is given by 
    \[
      \eta(C_\lambda)(\mu)=
      k_\lambda
      Q_\lambda^*(-\mu),
    \]
    where
    \[
    k_\lambda\defi (-1)^{|\lambda|}\lambda_1!\dotsm \lambda_{\ell(\lambda)}!
    \prod_{1\sle i<j\sle \ell(\lambda)}
    \frac{\lambda_i+\lambda_j}{\lambda_i-\lambda_j}
    .
    \]
  \end{Th}

  \noindent
  We would like to mention that the polynomials $Q_\lambda$ occur in a further different scenario related to the Lie superalgebra $\ger q(n)$. In \cite{sergeevv}*{Theorem 1.7}, Sergeev showed that the radial parts of the bispherical matrix coefficients on $\ger q(n)\times\ger q(n)$ with respect to the diagonal and twisted-diagonal embeddings of $\ger q(n)$ in $\ger q(n)\times\ger q(n)$, are Schur $Q$-polynomials. It will be interesting to explore possible connections between our work and Sergeev's result.

  We remark that it is possible to extend the results of the present paper and of Ref.~\cite{sahi-salmasian} to the common setting of multiplicity-free actions on Jordan superalgebras. This will be established in a forthcoming paper \cite{SahiSalSer}. In addition, recently Sahi and Salmasian~\cite{sahi-salmasian-ENS} constructed quadratic analogues of Capelli operators on Grassmannian manifolds by lifting the Capelli basis of~\cite{s94} via a double fibration. In the near future they plan to consider the analogous problem in the super setting.

  We conclude this introduction with a brief synopsis of our paper. In Section \ref{s:lsa}, we realise the Lie superalgebras relevant to us in terms of supermatrices. In Section \ref{s:evpol}, we identify the action of $\ger l$ on $\mathscr P(V)$, construct the Capelli basis, and determine the eigenvalue polynomials (\thmref{Th}{evpol-q}). Finally, in Section \ref{s:sphpol}, we study the open orbits in $V$ and $V^*$, show the existence of $\ger m$-invariant functionals for the simple summands of $\mathscr P(V)$, and prove that the spherical polynomials thus defined are the classical Schur $Q$-functions (\thmref{Th}{sphpol-q-hompart}).

  \medskip
  \noindent\emph{Acknowledgements.} 
  Alexander Alldridge wishes to thank the University of Ottawa, the Institute for Theoretical Physics at the University of Cologne, and the Department of Mathematics of the University of California at Berkeley for their hospitality during the preparation of this article. 

  The authors thank Vera Serganova and Weiqiang Wang for stimulating and fruitful conversations during the workshop, which paved the way for the present article. We also thank Alexander Sergeev for bringing Ref.~\cite{sergeevv} to our attention. 

  \section{Lie superalgebras}\label{s:lsa}

  The triple $(\ger l,\ger m,V)$ given in Equation \eqref{qqq} can be embedded inside the Lie superalgebra $\ger q(2n)$, which can be further embedded inside $\ger{gl}(2n|2n)$. This provides a concrete realisation which allows us to express the Lie superalgebras of interest as $4n\times 4n$ matrices. In order to describe it, it will be convenient to consider three commuting involutions of the algebra $\ger{gl}(2n|2n)$. To this end, first we equip the space $\mathrm{Mat}_{4n\times 4n}(\cplxs)$ of $4n\times 4n$ complex matrices with a Lie superalgebra structure isomorphic to $\ger{gl}(2n|2n)$. Instead of supermatrices in standard format, we prefer to consider those of the shape
  \begin{equation}\label{eq:gl-nonstd}
    x=
    \begin{Matrix}1
      A&B\\ C&D
    \end{Matrix},\quad A,B,C,D\in\ger{gl}(n|n).
  \end{equation}
  Equipped with the signed matrix commutator, the space of such matrices forms a Lie superalgebra $\ger g$ isomorphic to $\ger{gl}(2n|2n)$, the isomorphism being given by conjugation by the $4n\times 4n$ matrix
  \[
    \widetilde I\defi\begin{Matrix}1
      I&0&0&0\\0&0&I&0\\0&I&0&0\\0&0&0&I
    \end{Matrix},
  \] 
  where $I\defi I_{n\times n}$ denotes the $n\times n$ identity matrix. Next let $\sigma, \varphi$, and $\theta$  be three involutions on $\ger g$, given respectively by conjugation by the matrices
  \[
    \Sigma\defi\begin{Matrix}1
    I & 0 & 0 & 0\\
    0 & I & 0 & 0\\
    0 & 0 & -I & 0\\
    0 & 0 & 0 & -I
    \end{Matrix}
    \ ,\ 
    \Phi\defi\begin{Matrix}1
    0 & I & 0 & 0\\
    I & 0 & 0 & 0\\
    0 & 0 & 0 & I\\
    0 & 0 & I & 0
    \end{Matrix}
    \ ,\ 
    \Theta\defi\begin{Matrix}1
    0 & 0 & I & 0\\
    0 & 0 & 0 & I\\
    I & 0 & 0 & 0\\
    0 & I & 0 & 0
    \end{Matrix}.
  \]
  It is straightforward to verify that $\Sigma\Phi=\Phi\Sigma$, $\Phi\Theta=\Theta\Phi$, and $\Sigma\Theta=-\Theta\Sigma$. Hence, the involutions $\sigma,\vphi$, and $\theta$ commute with each other.

  \subsection{\texorpdfstring{The involution $\sigma$}{The involution sigma}}

  The subspace $\ger g^\sigma$ of fixed points  of $\sigma$ equals $\ger{gl}(n|n)\times\ger{gl}(n|n)$ with elements of the form 
  \[
    \begin{Matrix}1
      A&0\\0&D
    \end{Matrix}\text{ where }A,D\in\ger{gl}(n|n).
  \]
  The subspace $\ger g^{-\sigma}$ of fixed points of $-\sigma$ consists of matrices of the form 
  \[
    \begin{Matrix}1
      0&B\\C&0
    \end{Matrix}\text{ where }B,C\in\ger{gl}(n|n).
  \]
  Thus, as a super-vector space, $\ger g^{-\sigma}=\ger g^{-\sigma}_+\oplus\ger g^{-\sigma}_-$ where $\ger g^{-\sigma}_\pm$ are respectively  the spaces of $4n\times 4n$ matrices of the form
  \[
    \begin{Matrix}1
      0&B\\0&0
    \end{Matrix}\text{ and }
    \begin{Matrix}1
      0&0\\
      C&0
    \end{Matrix}
  \]
  with $B,C\in\ger{gl}(n|n)$. In fact, $\ger g^{-\sigma}_\pm$ are the $\pm 2$-eigenspaces of $\mathrm{ad}(\Sigma)$, where we think of $\Sigma$ as an element of $\ger g$. Therefore, $\ger g^{-\sigma}_\pm$ are abelian subalgebas of $\ger g$, and together with $\ger g^\sigma$, they form a $\ints$-grading of $\ger g$. The action of $\ger g^\sigma$ on $\ger g^{-\sigma}_+$ is given explicitly by
  \begin{equation}\label{eq:fpalg-action}
    \Bracks1{x,v}=
    \begin{Matrix}1
      0&AB-(-1)^{\Abs0B\Abs0D}BD\\0&0
    \end{Matrix},
  \end{equation}
  \fa homogeneous 
  \[
    x=
    \begin{Matrix}1
      A&0\\0&D
    \end{Matrix}\in\ger g^\sigma,\quad
    v=
    \begin{Matrix}1
      0&B\\0&0
    \end{Matrix}\in\ger g^{-\sigma}_+.
  \]
  In what follows, we set $U\defi\ger g^{-\sigma}_+$ and identify it as a super-vector space with $\ger{gl}(n|n)$, \via the map
  \begin{equation}\label{eq:u-id}
    \ger{gl}(n|n)\longrightarrow U:B\longmapsto
    \begin{Matrix}1
      0&B\\0&0
    \end{Matrix}.
  \end{equation}

  \subsection{\texorpdfstring{The involution $\vphi$}{The involution phi}}

  The involution $\varphi$ is induced by the parity reversing automorphism of $\cplxs^{2n|2n}$, and therefore the subalgebra $\ger g^\varphi$ of fixed points of $\varphi$  is isomorphic to  $\ger q(2n)$. It consists of all $x\in\ger g$ as in Equation \eqref{eq:gl-nonstd} \scth the blocks $A,B,C,D\in\ger q(n)$. Furthermore, 
  \[
    \ger l\defi \ger g^{\vphi,\sigma}\cong\ger q(n)\times\ger q(n).
  \]

  We define $V\defi\ger g^\vphi\cap U$. Then it is clear that 
  \[
    [\ger l,V]\subseteq\ger g^\vphi\cap[\ger g^\sigma,\ger g^{-\sigma}_+]\subseteq\ger g^\vphi\cap\ger g^{-\sigma}_+=V.
  \]
  The restriction of the map defined in Equation \eqref{eq:u-id} yields an identification of $V$ with a subspace of $\ger{gl}(n|n)$ which carries the structure of $\ger q(n)$.

  \subsection{\texorpdfstring{The involution $\theta$}{The involution theta}} 

  The algebra $\ger g^{\sigma,\theta}$ of fixed points of both $\sigma$ and $\theta$ is isomorphic to $\ger{gl}(n|n)$, and realised by supermatrices 
  \begin{equation}\label{eq:tau-fpalg}
    \begin{Matrix}1
      A&0\\0&A
    \end{Matrix},\quad A\in\ger{gl}(n|n).
  \end{equation}
  From Equation \eqref{eq:fpalg-action}, it follows that the action of $\ger g^{\sigma,\theta}$ on $U$ is precisely the adjoint action of $\ger{gl}(n|n)$. Furthermore,
  \[
     \ger m\defi\ger l^\theta=\ger g^{\vphi,\sigma,\theta}\cong \ger q(n).
  \] 
  It is realised by supermatrices of the form as in Equation \eqref{eq:tau-fpalg} where in addition $A\in\ger q(n)$. Moreover, the action of $\ger m$ on $V$ is precisely the adjoint action of $\ger q(n)$. 
  For the following lemma, let 
$e\in V\subseteq U$ be the element correponding to the matrix
\[
\begin{Matrix}1
0 & I_{2n\times 2n}\\ 0 & 0
\end{Matrix},
\]  
and  
  set 
  \[
    \ger s\defi\ger l^{-\theta}=\ger g^{\vphi,\sigma,-\theta}.
  \]

  \begin{Lem}[m-isotropy]
    Let $\mathsf p:\ger g^\sigma\longrightarrow U$ be the linear map defined by $x\longmapsto x\cdot e$.
    \begin{enumerate}[wide]
      \item\label{item:m-isotropy-i} The map $\mathsf p$ is a surjection onto $U$ with  kernel $\ger g^{\sigma,\theta}$. Its restriction to $\ger g^{\sigma,-\theta}$ is a linear isomorphism.
      \item\label{item:m-isotropy-ii} The restriction of $\mathsf p$ to  $\ger l$ is a surjection onto $V$ with kernel $\ger m$. Its restriction to  $\ger s$ is a linear isomorphism onto $V$.
      \item\label{item:m-isotropy-iii} Equipped with the binary operation $(x\cdot e)\circ (y\cdot e)\defi[x,[y,e]]$ for $x,y\in \ger g^{\sigma,-\theta}$ (respectively, $x,y\in\ger s$)
      the super-vector space $U$ (respectively, $V$) becomes a Jordan superalgebra.
    \end{enumerate}
  \end{Lem}

  \begin{proof}
    Parts \eqref{item:m-isotropy-i} and \eqref{item:m-isotropy-ii} follow from straightforward calculations. For \eqref{item:m-isotropy-iii}, it suffices to prove the statement for $U$. Note that according to Equation \eqref{eq:fpalg-action}, we have 
    \[
      x\cdot e=
      \begin{Matrix}1
        0&2A\\0&0
      \end{Matrix},
      \quad\forall
      x=
      \begin{Matrix}1
        A&0\\0&-A
      \end{Matrix}\in\ger g^{\sigma,-\theta}.
    \]
    Thus, a direct calculation shows that up to normalization, $A\circ B$ coincides for $A,B\in U$ with the super-anticommutator $AB+(-1)^{|A|\cdot|B|}BA$. 
  \end{proof}
 
  \section{The eigenvalue polynomials}\label{s:evpol}

  \subsection{\texorpdfstring{The action of $\ger l$ on polynomials}{The action of l on polynomials}} \label{Sec3.1}
  
  Let $\mathscr P(V)$ denote the superalgebra of super-polynomials on $V$. Recall that $\mathscr P(V)$ is by definition equal to $S(V^*)$. As $\ger g^\sigma$ acts on $V$, we obtain an induced locally finite $\ger g^\sigma$-action on $\mathscr P(V^*)=S(V)$. Similar statements apply to $\mathscr P(V)$.
  We shall identify this action in terms of differential operators. To that end, we consider the complex supermanifold $\aff(V)$, defined as the locally ringed space with underlying topological space $V_\ev$ and sheaf of superfunctions $\sh O_{\aff(V)}\defi\sh H_{V_\ev}\otimes\bigwedge(V_\odd)^*$, where $\sh H$ denotes the sheaf of holomorphic functions. There is a natural inclusion
  \[
    V^*\longrightarrow\Gamma(\sh O_{\aff(V)})
  \]
  ($\Gamma$ denoting global sections), allowing us to identify linear forms on $V$ with certain superfunctions on $\aff(V)$. In particular, $\mathscr P(V)$ is a subsuperalgebra of $\Gamma(\sh O_{\aff(V)})$.

  Recall that on a supermanifold $X$, the vector fields on $X$, defined on an open set $O\subseteq X_0$ of the underlying topological space, are defined to be the superderivations of $\sh O_X|_O$ \cite{ahw-orbits}*{Definition 4.1}. For any homogeneous basis $(x_a)$ of $V$, the dual basis $(x^a)$ is a coordinate system on $\aff(V)$, and \cite{ahw-orbits}*{Proposition 4.5} there are unique (and globally defined) vector fields $\frac\partial{\partial x^a}$ on $\aff(V)$ of parity $\Abs0{x^a}$, determined by
  \[
    \frac\partial{\partial x^a}(x^b)=\delta_{ab},\quad\forall a,b.
  \]
  The linear action of $\ger l$ on $V$ determines, for $v\in\ger l$, vector fields $\mathbf a_v$ on $\aff(V)$ by
  \begin{equation}\label{eq:fund-vf}
    \mathbf a_v(x^b)=-v\cdot x^b,\quad\forall b.
  \end{equation}
  The sign stems from the fact that these are the fundamental vector fields for a Lie supergroup action on $\aff(V)$, as we shall see later. By construction, for any $v\in\ger l$, the action of $-a_v$ on $\mathscr P(V)$ coincides with the action of $v$ defined in the first paragraph of this subsection. 

  We now make this action explicit. Let $(e_k,\eps_k)_{k=1,\dotsc,n}$ be the standard basis of $\cplxs^{n|n}$. A homogeneous basis of the super-vector space $V\cong\ger q(n)$ is determined by
  \[
    u_{k\ell}\defi
    \begin{Matrix}1
      E_{k\ell}&0\\0&E_{k\ell}
    \end{Matrix},\quad
    \xi_{k\ell}\defi
    \begin{Matrix}1
      0&E_{k\ell}\\E_{k\ell}&0
    \end{Matrix},
  \]
  where $k,\ell=1,\dotsc,n$ and $E_{k\ell}$ are the usual elementary matrices. Then 
  \begin{align*}
    u_{k\ell}(e_j)&=\delta_{\ell j}e_k, &u_{k\ell}(\eps_j)&=\delta_{\ell j}\eps_k,\\
    \xi_{k\ell}(e_j)&=\delta_{\ell j}\eps_k, &\xi_{k\ell}(\eps_j)&=\delta_{\ell j}e_k.
  \end{align*}
  Let $(u^{k\ell},\xi^{k\ell})$ be the dual basis of $\ger q(n)^*=V^*$. This determines vector fields 
  \[
    \frac\partial{\partial u^{k\ell}},\quad \frac\partial{\partial \xi^{k\ell}}
  \]
  on $\aff(V)$, by the recipe given above. 

  Moreover, let $(a_{k\ell},\alpha_{k\ell})$ and $(b_{k\ell},\beta_{k\ell})$, respectively, be the copies of $(u_{kl},\xi_{k\ell})$ in the first and second factor of $\ger l=\ger q(n)\times\ger q(n)$. By Equation \eqref{eq:fpalg-action}, as a module over the second factor, $V$ is isomorphic to $\cplxs^n\otimes(\cplxs^{n|n})^*$ where $\ger q(n)$ acts on $\cplxs^{n|n}$ in the standard way. Hence, the second factor acts on $V^*$ as on $\cplxs^n\otimes\cplxs^{n|n}$:
  \begin{align*}
    b_{k\ell}(u^{pq})&=\delta_{\ell q}u^{pk}&b_{k\ell}(\xi^{pq})&=\delta_{\ell q}\xi^{pk},\\
    \beta_{k\ell}(u^{pq})&=\delta_{\ell q}\xi^{pk}&\beta_{k\ell}(\xi^{pq})&=\delta_{\ell q}u^{pk}.
  \end{align*}
  It follows that 
  \begin{equation}\label{eq:secondfactor-action}
    -\!{\mathbf a_{b^{k\ell}}}=\sum_{p=1}^nu^{pk}\frac\partial{\partial u^{p\ell}}+\xi^{pk}\frac\partial{\partial\xi^{p\ell}},\quad
    -{\mathbf a_{\beta^{k\ell}}}=\sum_{p=1}^n\xi^{pk}\frac\partial{\partial u^{p\ell}}+u^{pk}\frac\partial{\partial\xi^{p\ell}}.
  \end{equation}

  Similarly, the first factor of $\ger l$ acts on $V$ as on $\cplxs^{n|n}\otimes\cplxs^n$, and hence on $V^*$ as on $(\cplxs^{n|n})^*\otimes\cplxs^n$. Reasoning as above, this implies
  \begin{equation}\label{eq:firstfactor-action}
    -\!{\mathbf a_{a^{k\ell}}}=\sum_{p=1}^nu^{p\ell}\frac\partial{\partial u^{pk}}+\xi^{p\ell}\frac\partial{\partial\xi^{pk}},\quad
    -{\mathbf a_{\alpha^{k\ell}}}=\sum_{p=1}^n-\xi^{p\ell}\frac\partial{\partial u^{pk}}+u^{p\ell}\frac\partial{\partial\xi^{pk}}.
  \end{equation}

  We will presently decompose the $\ger l$-module $\mathscr P(V)$. To that end, we introduce a labelling set. Let $\Lambda$ be the set of \Define{partitions}, that is, of all finite sequences $\lambda=(\lambda_1,\dotsc,\lambda_m)$ of non-negative integers $\lambda_j$ \scth $\lambda_1\sge\lambda_2\sge\dotsm\sge\lambda_m$. Here, we identify $\lambda$ with any partition $(\lambda_1,\dotsc,\lambda_m,0,\dotsc,0)$ obtained from $\lambda$ by appending a finite number of zeros at its tail. If $\lambda$ can be written in the form $(\lambda_1,\dotsc,\lambda_\ell)$ where $\lambda_\ell>0$, then we say $\lambda$ has \Define{length} $\ell(\lambda)=\ell$. Let $\Lambda^\ell\subseteq\Lambda$ be the set of partitions of length $\ell(\lambda)\sle \ell$. We also set $\Abs0\lambda\defi\lambda_1+\dotsm+\lambda_\ell$ if $\ell(\lambda)\sle \ell$. A partition $\lambda$ of length $\ell(\lambda)=\ell$ is called \Define{strict} if $\lambda_1>\dotsm>\lambda_\ell>0$. The set of all strict partitions will be denoted by $\Lambda_{>0}$, and we write $\Lambda^\ell_{>0}$ for the set of strict partitions of length at most $\ell$; that is, $\Lambda^\ell_{>0}\defi\Lambda^\ell\cap\Lambda_{>0}$.

  For every strict partition $\lambda$ such that $\ell(\lambda)\sle n$, let $F_\lambda$ be the $\ger q(n)$-highest weight module with highest weight $\lambda_1\eps_1+\dotsm+\lambda_n\eps_n$, where the $\eps_i$ are the standard characters of the even part of the Cartan subalgebra of $\ger q(n)$. For every strict partition $\lambda$, set $\delta(\lambda)\defi0$ when $\ell(\lambda)$ is even, and $\delta(\lambda)\defi1$ otherwise. We now define an $\ger l$-module $V_\lambda$ as follows. It is shown in \cite{cw}*{Section 2} that, as an $\ger l$-module, the exterior tensor product  $(F_\lambda)^*\boxtimes F_\lambda$ is irreducible when $\delta(\lambda)=0$, and decomposes into a direct sum of two irreducible  isomorphic $\ger l$-modules (via an odd map) if  $\delta(\lambda)=1$. Following the notation of Ref.~\cite{cw}, we set
  \[
    V_\lambda\defi\frac1{2^{\delta(\lambda)}}\Parens1{F_\lambda)^*\boxtimes F_\lambda},
  \]
  that is, we take $V_\lambda$ to be the irreducible component of $(F_\lambda)^*\boxtimes F_\lambda$ that appears in the decomposition of the super-polynomial algebra over the natural  $(\ger q(n),\ger q(n))$-module (see \thmref{Prop}{polyact-ssmf}). The $\ger l$-module $V_\lambda$ is always of type $\mathsf M$, that is, it is irreducible as an ungraded representation. It follows that in the $\mathbb Z_2$-graded sense, 
  \begin{equation}\label{hmVlVl1}
    \GHom[_{\ger l}]1{V_\lambda,V_\mu}=\delta_{\lambda\mu}\cdot\cplxs,
  \end{equation}
  where \emph{a priori}, $\underline{\mathrm{Hom}}_{\ger l}$ denotes the set of \emph{all} $\ger l$-equivariant linear maps (of any parity). In particular, all non-zero $\ger l$-equivariant endomorphisms of $V_\lambda$ are even.
  
  \begin{Prop}[polyact-ssmf]
    Under the action of $\ger l$, $\mathscr P^k(V)$ decomposes as the multiplicity-free direct sum of simple modules $V_\lambda$, where $\lambda$ ranges over elements of $\Lambda^n_{>0}$ which satisfy $\Abs0\lambda=k$.  
  \end{Prop}

  \begin{proof}
    Recall that $\mathscr P(V)\cong S(V^*)$ as $\ger l$-modules. The proposition follows from the description of the actions of the left and right copies of $\ger q(n)$ on $V$ given above, and the results of \cite{cw}*{Section 3}.
  \end{proof}

  \subsection{Invariant polynomial differential operators}
  \label{Sec3.2}

  On a complex supermanifold $X$, the \Define{differential operators} on $X$ defined on an open set $O\subseteq X_0$ of the underlying topological space are generated as a subsuperalgebra of the $\cplxs$-linear endomorphisms of $\sh O_X|_O$ by vector fields and functions. This gives a $\cplxs$-algebra sheaf $\mathscr D_X$ that is an $\sh O_X$-bisupermodule, filtered by order. 

  Here, a differential operator is of \Define{order} $\sle n$ if it can be expressed as a product of some functions and at most $n$ vector fields. A differential operator of order $\sle n$ is uniquely determined by its action on monomials of order $\sle n$ in some given system of coordinate functions. This follows in the usual way from the Hadamard Lemma \cite{leites}*{Lemma 2.1.8} and implies that $\mathscr D_X$ is locally free as a left $\sh O_X$-supermodule.

  For the supermanifold $\aff(V)$, we have a $\cplxs$-superalgebra map 
  \[
    \partial:S(V)\longrightarrow\Gamma(\mathscr D_{\aff(V)}).
  \]
  It is determined by the linear map which sends any homogeneous $v\in V$ to the unique vector field $\partial(v)$ \scth
  \begin{equation}\label{eq:del-def}
    \partial(v)(\mu)=(-1)^{\Abs0\mu\Abs0v}\mu(v)
  \end{equation}
  \fa homogeneous $\mu\in V^*$. If $(x_a)$ is a homogeneous basis of $V$, then
  \[
    \partial(x_a)=(-1)^{|x_a|}\frac\partial{\partial x^a}.
  \]
  The image of $\partial$, denoted by $D(V)$, is the superalgebra of \emph{constant-coefficient} differential operators on $V$. The map
  \begin{equation}\label{eq:pd-iso}
    \mathscr P(V)\otimes D(V)\longrightarrow\Gamma(\mathscr  D_{\aff(V)}):p\otimes D\longmapsto pD
  \end{equation}
  is an isomorphism onto a $\mathscr P(V)$-submodule of $\Gamma(\mathscr D_{\aff(V)})$ denoted by $\mathscr{PD}(V)$. Indeed, $\mathscr{PD}(V)$ is a subsuperalgebra, 
  the algebra of \Define{polynomial differential operators}. As $\ger l$ acts by linear vector fields, we have the bracket relation 
  \[
    [-\mathbf a_x,p]=-\mathbf a_x(p)=x\cdot p,\quad\forall x\in\ger l,p\in\mathscr P(V).
  \]
  The bilinear form 
    \begin{equation}
    \label{poaringn}
      \mathscr P^k(V)\otimes D^k(V)\longrightarrow\cplxs:p\otimes D\longmapsto Dp
    \end{equation}
  is a non-degenerate pairing which is $\ger l$-equivariant, and therefore results in a canonical $\ger l$-module isomorphism $D^k(V)\cong \mathscr P^k(V)^*$. Hence the following corollary to \thmref{Prop}{polyact-ssmf} holds.

  \begin{Cor}[{cor3.2}]
    The space $\mathscr{PD}(V)^\ger l$ of $\ger l$-invariant polynomial differential operators decomposes as follows:
    \[
      \mathscr{PD}(V)^\ger l=\bigoplus_{\lambda\in\Lambda^n_{>0}}
      \Parens1{V_\lambda\otimes (V_\lambda)^*}^\ger l.
    \]
  \end{Cor}

  There is a natural $\ger l$-equivariant isomorphism $V_\lambda\otimes (V_\lambda)^*\cong\GEnd[_\cplxs]0{V_\lambda}$, so the identity element ${\id}_{V_\lambda}$ of $V_\lambda$ determines an even $\ger l$-invariant polynomial differential operator
  \[
    D_\lambda\in\mathscr{PD}(V)^\ger l.
  \]
  That is, if $(p_j)$ is a homogeneous basis of $V_\lambda\subseteq\mathscr P(V)$ and $(D_j)$ is its dual basis for $(V_\lambda)^*\subseteq D(V)$, then 
  \begin{equation}\label{eq:dlambda-def}
    D_\lambda\defi\sum_j p_jD_j.
  \end{equation}
  Schur's Lemma implies that for every $\mu\in\Lambda^n_{>0}$, there is a complex scalar $c_\lambda(\mu)\in\cplxs$ \scth $D_\lambda$ acts by the scalar $c_\lambda(\mu)$ on $V_\lambda\subseteq\mathscr P(V)$. 

  \begin{Cor}[ev-props]
    The operators $D_\lambda$, where $\lambda$ ranges over $\Lambda^n_{>0}$, form a basis of the space $\mathscr{PD}(V)^\ger l$ of $\ger l$-invariant differential operators. Moreover, $c_\lambda(\lambda)=1$, while $c_\lambda(\mu)=0$ whenever $\Abs0\mu\sle\Abs0\lambda$ and $\mu\neq\lambda$.
  \end{Cor}

  \begin{proof}
    The first statement follows from \thmref{Cor}{cor3.2} and Equation \eqref{hmVlVl1}. Since the order of the operator $D_\lambda$ is not $\sle\Abs0\lambda-1$, it vanishes on $\mathscr P^k(V)$ for $k<\Abs0\lambda$. Next assume $k=\Abs0\lambda=\Abs0\mu$. If $\lambda\neq \mu$ and $D_\lambda$ does not vanish on $V_\mu$, then the restriction of the bilinear form \eqref{poaringn} to $V_\mu\times (V_\lambda)^*$ will be a nonzero $\ger l$-equivariant form, hence $V_\lambda\cong V_\mu$, which is a contradiction.

    It remains to compute the action of $D_\lambda$ on $V_\lambda$. Let $(p_j)$ and $D_j$ denote the dual bases of $V_\lambda$ and $(V_\lambda^*)$. Then $D_\lambda p_i=\sum_jp_jD_jp_i=\sum_jp_j\delta_{ij}=p_i$ for every $i$.
  \end{proof}

  To determine $c_\lambda(\mu)$ (which will be done in the next subsection), we first need to see that it extends to a $Q$-symmetric polynomial. To that end, let $\Gamma_n$ denote the ring of \Define{$Q$-symmetric polynomials}, that is, $n$-variable symmetric polynomials $p(x_1,\dotsc,x_n)$ such that $p(t,-t,x_3,\dotsc,x_n)$ does not depend on $t$. (For $n=1$, the latter condition is vacuous.) 

  \begin{Prop}[ev-pol]
    \Fa $\lambda\in\Lambda^n_{>0}$, there exists a polynomial $q_\lambda^*\in\Gamma_n$ of degree at most $\Abs0\lambda$ \scth $c_\lambda(\mu)=q_\lambda^*(\mu)$ \fa $\mu\in\Lambda^n_{>0}$.
  \end{Prop}

  \begin{proof}
    Recall that $L$ and $R$, respectively, denote the actions of the first and second factor of $\ger l=\ger q(n)\times\ger q(n)$ on $\mathscr P(V)$. As $\mathscr P(V)$ is the multiplicity-free direct sum of simple modules of the universal enveloping algebra $\ger U(\ger l)=\ger U(\ger q(n))\otimes\ger U(\ger q(n))$, it follows that $L(\ger U(\ger q(n)))$ and  $R(\ger U(\ger q(n)))$ are mutual commutants in $\mathscr{PD}(V)$ (this double commutant property is also mentioned in \cite{cw}). In particular, we have 
    \[
      \mathscr{PD}(V)^\ger l\subseteq R(\ger U(\ger q(n))).
    \]
    As $R$ is faithful, it follows that in fact 
    \[
      \mathscr{PD}(V)^\ger l=R(\sh Z(\ger q(n))),
    \]
    where $\sh Z$ denotes the centre of $\ger U$. The latter statement also follows from the explicit construction of the Capelli operators in the work of Nazarov \cite{nazarov}.

    Furthermore, the simple module $V_\lambda$ occuring in the decomposition of $\mathscr P(V)$ is contained in the external tensor product $(F_\lambda)^*\boxtimes F_\lambda$. By Sergeev's Harish-Chandra isomorphism for $\ger q(n)$ \cite{sergeev}*{Theorem 3} (see also \cite{cwbook}*{Theorem 2.46}), for any $u\in\sh Z(\ger q(n))$, there exists a $Q$-symmetric polynomial $q_u^*\in\Gamma_n$ \scth \fa $\lambda\in\Lambda^n_{>0}$, $u$ acts on $F_\lambda$ by $q^*_u(\lambda)$.

    Fix $\lambda\in\Lambda^n_{>0}$. Then the order of $D_\lambda$ is $\sle\Abs0\lambda$. As $L$ is faithful, there is a unique $z_\lambda\in\sh Z(\ger q(n))$ \scth $L(z_\lambda)=D_\lambda$. Since $\ger q(n)\subseteq\ger l$ acts by linear vector fields, $z_\lambda$ lies in the $\Abs0\lambda$-th part of the standard increasing filtration of $\ger U(\ger q(n))$. Then $q_\lambda^*\defi q^*_{z_\lambda}$ has degree at most $\Abs0\lambda$, see \cites{sergeev,cwbook}. The assertion follows.
  \end{proof}

  \begin{Def}
    We call $q_\lambda^*$ the \Define{eigenvalue polynomial} of $D_\lambda$ for $\lambda\in\Lambda^n_{>0}$.
  \end{Def}

  \subsection{\texorpdfstring{Schur $Q$-functions}{Schur Q-functions}}

  Our next goal is to identify the eigenvalue polynomials $q_\lambda^*\in\Gamma_n$. We first recall the definitions of certain elements $Q_\lambda$ of $\Gamma_n$, called the \emph{$Q$-functions} of Schur, and their shifted analogues, the \Define{factorial $Q$-functions} $Q^*_\lambda$ originally defined by Okounkov, see Ref.~\cite{ivanov97}.

  Given a sequence $(a_n)_{n\sge1}$ of complex numbers, we define for any non-negative integer $k$ the \Define{$k$th generalized power} of $x$ by
  \[
    (x\mid a)^k\defi\prod_{i=1}^k (x-a_i),
  \]
  where we set $(x|a)^0=1$.
  For every $\lambda\in\Lambda^n_{>0}$, we set
  \begin{align*}
    F_\lambda(x_1,\dotsc,x_n\mid a)&\defi\prod_{i=1}^{\ell(\lambda)}(x\mid a)^{\lambda_i}\prod_{\substack{i\sle\ell(\lambda)\\ i<j\sle n}}\frac{x_i+x_j}{x_i-x_j}.
  \end{align*}
  We now define
  \begin{align*}
    Q_\lambda(x_1,\dotsc,x_n\mid a)
    &
    \defi \frac{2^{\ell(\lambda)}}{(n-\ell(\lambda))!}
    \sum_{\sigma\in S_n}
    F_\lambda(x_{\sigma(1)},\dots ,x_{\sigma(n)}\mid a).
  \end{align*}
  We remark that $Q_\lambda(x_1,\dotsc,x_n\mid  a)$ can be expressed as a ratio of an antisymmetric polynomial by the Vandermonde polynomial, and therefore it is also a polynomial. Two special cases of interest are
  \[
    Q_\lambda(x_1,\dotsc,x_n)\defi Q_\lambda(x_1,\dotsc,x_n\mid \underline{0}),
  \] 
  where $\underline{0}\defi(0,0,0,\dotsc)$, and 
  \[
    Q_\lambda^*(x_1,\dotsc,x_n)\defi Q_\lambda(x_1,\dotsc,x_n\mid \delta),
  \]
  where $\delta\defi(0,1,2,3,\dotsc)$, called, respectively, the \Define{Schur $Q$-function}, and the \Define{factorial Schur $Q$-function}. 

  \begin{Prop}[PrptiQlQlstar]
    Let $\lambda\in\Lambda^n_{>0}$. 
    \begin{enumerate}[wide]
      \item\label{item:PrptiQlQlstar-i} We have $Q_\lambda,Q_\lambda^*\in\Gamma_n$. Furthermore, $Q_\lambda$ is homogeneous of degree $\Abs0\lambda$ and  
      \[
        \deg(Q_\lambda^*-Q_\lambda)<\Abs0\lambda.
      \]
      \item\label{item:PrptiQlQlstar-ii} Both $(Q_\mu)_{\mu\in\Lambda^n_{>0}}$ and $(Q^*_\mu)_{\mu\in\Lambda^n_{>0}}$ are bases for the vector space $\Gamma_n$. 
      \item\label{item:PrptiQlQlstar-iii} We have $Q_\lambda^*(\mu)=0$ for every $\mu\in\Lambda^n_{>0}$ \scth $\Abs0\mu\sle\Abs0\lambda$, $\lambda\neq\mu$; moreover, 
      \[
        Q_\lambda^*(\lambda)=Q_\lambda(\lambda)=H(\lambda)\defi\lambda!\prod_{i<j}\frac{\lambda_i+\lambda_j}{\lambda_i-\lambda_j},
      \] 
      where $\lambda!\defi\lambda_1!\dotsm\lambda_n!$, and it is understood that the product extends only up to $\ell(\lambda)$.
      \item\label{item:PrptiQlQlstar-iv} The unique element of $\Gamma_n$ of degree at most $\Abs0\lambda$ which satisfies Equation \eqref{item:PrptiQlQlstar-iii} is precisely $Q_\lambda^*$.
    \end{enumerate}
  \end{Prop}

  \begin{proof}
    Parts \eqref{item:PrptiQlQlstar-i}--\eqref{item:PrptiQlQlstar-iii} follow from \cite{ivanov97}*{\S~1} (see also \cite{Stembridge}*{\S~6}). Thus, we only sketch the argument for item \eqref{item:PrptiQlQlstar-iv}. Set $k\defi\Abs0\lambda$ and let $\Gamma_n(k)\subseteq\Gamma_n$ consist of the polynomials of degree at most $k$. From item \eqref{item:PrptiQlQlstar-ii}, it follows that $\dim(\Gamma_n(k))=\#\Lambda^n_{>0}(k)$, where
    \[
      \Lambda^n_{>0}(k)\defi\Set1{\lambda\in\Lambda^n_{>0}}{\Abs0\lambda\sle k}.
    \]
    For every $\mu\in\Lambda^n_{>0}(k)$, we consider the linear functional $\eps_\mu\in\Gamma_n^*$ defined by
    \[
      \eps_\mu(f)\defi f(\mu),\quad\forall f\in\Gamma_n.
    \]
    Let the order $\subseteq$ on $\Lambda$ be defined by 
    \[
      \mu\subseteq\nu\quad\text{if and only if}\quad\forall j:\mu_j\sle\nu_j.
    \]
    We choose a total order $\preceq$ on $\Lambda^n_{>0}(k)$ \scth $\mu\prec\nu$ if either $\Abs0\mu<\Abs0\nu$ or $\Abs0\mu=\Abs0\nu$ and $\nu\nsubseteq\mu$. Then by \cite{ivanov97}*{Proposition 1.16}, in terms of $\preceq$, the $k\times k$ matrix
    \[
      \Parens1{\eps_\mu(Q_\nu^*)}_{\mu,\nu\in\Lambda_{>0}^n(k)}
    \] 
    is triangular with no zeros on the diagonal, and therefore, invertible. It follows that the linear system
    \[
      \eps_\mu(f)=
      \begin{cases}
        H(\lambda)&\text{if }\mu=\lambda,\\
        0&\text{if }\mu\neq\lambda.
      \end{cases}
    \] 
    has a unique solution in $\Gamma_n(k)$. In view of item \eqref{item:PrptiQlQlstar-iii}, this solution is $Q_\lambda^*$.
  \end{proof}

  \begin{Rem}
    The polynomials $Q_\lambda$ also appear in Ref.~\cite{Macdonald}. They were initially introduced by Schur in connection with the projective representations of the symmetric group. According to \cite{ivanov97}*{Proposition 1.16}, the vanishing property of \thmref{Prop}{PrptiQlQlstar} \eqref{item:PrptiQlQlstar-iii} holds in a stronger form: $Q_\lambda^*(\mu)=0$ whenever $\lambda\nsubseteq\mu$. 
  \end{Rem}

  \begin{Th}[evpol-q]
    Let $q_\lambda^*$ for $\lambda\in\Lambda^n_{>0}$ be the eigenvalue polynomial of \thmref{Prop}{ev-pol}. Then 
    \[
      q_\lambda^*=\frac1{\lambda!}\prod_{i<j}\frac{\lambda_i-\lambda_j}{\lambda_i+\lambda_j}\,Q^*_\lambda.
    \]
  \end{Th}

  \begin{proof}
    By \thmref{Prop}{ev-pol},  $q_\lambda^*\in\Gamma_n$ satisfies the same degree, symmetry, and vanishing properties as $Q_\lambda^*$, see \thmref{Prop}{PrptiQlQlstar} \eqref{item:PrptiQlQlstar-iii}. The statement now follows from  \thmref{Prop}{PrptiQlQlstar} \eqref{item:PrptiQlQlstar-iv}.
  \end{proof}

  \begin{Rem}
    Nazarov \cite{nazarov}*{Proposition 4.8} constructs certain Capelli elements  $C_\lambda\in\mathcal Z(\ger q(n))$, and proves that the top degree part of their Harish-Chandra image  is the classical Schur $Q$-function. Even though he refers to Ivanov's work \cite{ivanov97} for the eigenvalue polynomials (see \cite{nazarov}*{p.~871}), his results do not include \thmref{Th}{evpol-q} explicitly.
  \end{Rem}

  \section{The spherical polynomials}\label{s:sphpol}

  In this final section, we define $\ger m$-spherical polynomials associated with the representations $V_\lambda$. We will show that up to a scalar multiple, these spherical polynomials are indeed the Schur $Q$-functions. As the Schur $Q$-function $Q_\lambda$ is the top-degree homogeneous part of $Q_\lambda^*$, our strategy is to prove that the spherical polynomials are equal to the top-degree homogeneous parts of the eigenvalue polynomials $q_\lambda^*$. 
  
  \subsection{Open orbits}

  We begin by globalising the action of $\ger l$ on $V$. Let $G$ be the complex Lie supergroup corresponding \via \cite{ccf}*{Theorem 7.4.5} to the supergroup pair $(\ger g,G_0)$, where $G_0$ consists of the matrices of the form
  \[
    \begin{Matrix}1
      a&0&b&0\\0&a'&0&b'\\c&0&d&0\\0&c'&0&d'
    \end{Matrix},\text{ such that } 
    \begin{Matrix}1
      a&b\\c&d
    \end{Matrix},
    \begin{Matrix}1
      a'&b'\\c'&d'
    \end{Matrix}\in\mathrm{GL}(2n,\cplxs),
  \]
  and acts on $\ger g$ by conjugation. Thus, $G_0\cong\GL(2n,\cplxs)\times\GL(2n,\cplxs)$. The automorphisms $\vphi$, $\sigma$, and $\theta$ of $\ger g$ integrate to $G$, as their restrictions to $\ger g_\ev$ integrate to $G_0$ and
  \[
    \Ad(\phi(g))(\phi(x))=\phi(\Ad(g))x,\quad\forall x\in\ger g,g\in G_0,
  \]
  for every choice of $\phi\in\{\vphi,\sigma,\theta\}$ \cite{ccf}*{Theorem 7.4.5}. In particular, there are fixed subsupergroups $G^\sigma$, $L\defi G^{\sigma,\vphi}$, and $M\defi G^{\sigma,\vphi,\theta}$ of $G$. One sees directly that 
  \begin{gather*}
    G_0^\sigma=\mathrm{GL}(n,\cplxs)\times\mathrm{GL}(n,\cplxs)\times\mathrm{GL}(n,\cplxs)\times\mathrm{GL}(n,\cplxs),\\
    L_0=\mathrm{GL}(n,\cplxs)\times\mathrm{GL}(n,\cplxs),\quad M_0=\mathrm{GL}(n,\cplxs)
  \end{gather*}
  are connected. 

  Write the elements of $G^\sigma_0$ in the form
  \[
    (A,D),\quad A=
    \begin{Matrix}1
      a&0\\0&a'
    \end{Matrix},\quad D=
    \begin{Matrix}1
      d&0\\0&d'
    \end{Matrix},\quad a,a',d,d'\in\mathrm{GL}(n,\cplxs).
  \]
  Then the action 
  \[
    G_0^\sigma\times U\longrightarrow U:(A,D,x)\longmapsto AxD^{-1}
  \]
  integrates the action of $\ger g_\ev^\sigma$ on $U=\ger{gl}(n|n)$ in such a way that 
  \[
    g\cdot (x\cdot v)=\Ad(g)(x)\cdot(g\cdot v),\quad\forall g\in G_0^\sigma,x\in\ger g^\sigma,u\in U,
  \]
  see Equation \eqref{eq:fpalg-action}. This gives rise to a Lie supergroup action
  \[
  G^\sigma\times \aff(U)\longrightarrow\aff(U),
  \]
  where $\aff(U)$ is the complex supermanifold corresponding to $U$, as in Section \ref{Sec3.1}.
  Passing to $\vphi$-fixed points, this implies that there is an action 
  \[
    \mathbf a:L\times\aff(V)\longrightarrow\aff(V)
  \]
  \scth the fundamental vector fields
  \[
    \mathbf a_x\defi(x\otimes 1)\circ \mathbf a^\sharp,\quad x\in\ger l,
  \]
  are those determined by Equation \eqref{eq:fund-vf}, see \cite{ccf}*{Proposition 8.3.3}. Here, $\mathbf a^\sharp$ is the  morphism of sheaves of superalgebras that corresponds to $\mathbf a$. We equally obtain an $L$-action on $\aff(V^*)$ integrating the given $\ger l$-action.

  Recall from Section \ref{Sec3.1} that $\aff(V)=\Parens1{V_\ev,\sh O_{\aff(V)}}$ as a locally ringed space.

  \begin{Prop}[openorb]
    The $L$-orbit $L\cdot e$ through the point $e\in V_\ev$ is the homogeneous supermanifold $L/M$. The canonical morphism $j_V:L\cdot e\longrightarrow \aff(V)$ is an $L$-equivariant open embedding. 
  \end{Prop}

  \begin{proof}
    As $M$ is connected, and by \thmref{Lem}{m-isotropy}, the isotropy of $L$ at $e$ is $M$, see \cite{ccf}*{Proposition 8.4.7}. Thus, $L\cdot e=L/M$ by definition, and the orbit morphism $\mathbf a_e:L\longrightarrow\aff(V)$ factors into $\pi:L\longrightarrow L/M$ and $j_V:L\cdot e=L/M\longrightarrow\aff(V)$ where $\pi$ is a surjective submersion and $j_V$ is an injective immersion, both $L$-equivariant, by \cite{ahw-orbits}*{Theorem 4.24}. By \thmref{Lem}{m-isotropy}, the kernel of 
    \[
      T_1(\mathbf a_e):T_1(L)=\ger l\longrightarrow T_e(\aff(V))=V
    \]
    is $\ger m$ and its restriction to $\ger s$ an isomorphism $\ger s\longrightarrow V$. Thus, the inverse function theorem \cite{leites}*{Theorem 2.3.1} applies, and $j_V$ is an open embedding. 
  \end{proof}

  \begin{Rem}
    It is not hard to identify the open orbit from \thmref{Prop}{openorb}. Indeed, identifying $V$ with $\ger q(n)$ \via the isomorphism from Equation \eqref{eq:u-id}, $L\cdot e$ equals $Q(n,\cplxs)\subseteq\aff(\ger q(n))$, defined to be the open subspace corresponding to the open subset $\GL(n,\cplxs)\subseteq\ger{gl}(n)=\ger q(n)_\ev$. Observe that $Q(n,\cplxs)$ has the structure of a complex Lie supergroup with Lie superalgebra $\ger q(n)$.

    The orbit morphism $L\longrightarrow L\cdot e$ is given by 
    \[
      (A,D)\longmapsto AD^{-1}
    \]
    on $T$-valued points, where $T$ is any complex supermanifold. 
  \end{Rem}

  \subsection{Spherical vectors and polynomials}\label{SpheP}

  Let $\lambda\in\Lambda^n_{>0}$. Define $p_\lambda^*\in(V_\lambda)^*$ by
  \begin{equation}\label{eq:sphpol-def}
    \Dual0{p_\lambda^*}\ell\defi (-1)^kk!\,j_V^\sharp(\ell)(e),\quad\forall\ell\in V_\lambda\subseteq\mathscr P(V),
  \end{equation}
  where $\Dual0\cdot\cdot$ is the standard pairing of $(V_\lambda)^*$ with $V_\lambda$, and $j_V^\sharp$ is pullback along the morphism $j_V:L\cdot e\longrightarrow\aff(V)$.

  \begin{Prop}[sphvects]
    For every $\lambda\in\Lambda^n_{>0}$, $p_\lambda^*$ is the up to scalars unique non-zero $\ger m$-invariant vector of $(V_\lambda)^*$. 
  \end{Prop}

  \begin{proof}
    Uniqueness of the $\ger m$-invariant in $V_\lambda$ follows from the fact that by the hom-tensor adjunction formula
    \[
      ((V_\lambda)^*)^\ger m=\GHom[_\ger m]0{V_\lambda,\cplxs}\cong\frac1{2^{\delta(\lambda)}}\,\GHom[_{\ger m}]0{F_\lambda,F_\lambda},
    \] 
    combined with the fact that $\dim\GEnd[_\ger m]0{F_\lambda}=2^{\delta(\lambda)}$. (Recall that $F_\lambda$ is of type $\mathsf Q$ if and only if $\delta(\lambda)=1$.)
    
    Next we prove that $p_\lambda^*$ is the desired $\ger m$-invariant. Since $j_V$ is $L$-equivariant and $M$ fixes $e$, it follows that $p_\lambda^*$ is $\ger m$-invariant. We need to see that $p_\lambda^*\neq0$. Above, we have noted the commutative diagram
    \[
      \begin{tikzcd}
        L\dar[swap]{\pi}\arrow{rd}{\mathbf a_e}\\
        L\cdot e\rar[swap]{j_V}&\aff(V)
      \end{tikzcd}
    \]
    where $j_V$ is an open embedding, $\pi$ is a surjective submersion, and 
    \[
      \mathbf a_e=\mathbf a\circ({\id}\times j_e),
    \]
    where $j_e$ is the embedding $*\longrightarrow\aff(V)$ of the singleton space defined by $e$. The action of $L$ is linear, \ie 
    \[
      \mathbf a^\sharp(V^*)\subseteq\Gamma(\sh O_L)\otimes V^*,
    \]
    and hence, because $\mathbf a^\sharp$ is a superalgebra morphism,  
    \[
      \mathbf a^\sharp(\mathscr P(V))\subseteq\Gamma(\sh O_L)\otimes\mathscr P(V).
    \]
    If $W\subseteq\mathscr P(V)$ is a graded subspace invariant under $L_0$ and $\ger l$, then 
    \[
      \mathbf a^\sharp(W)\subseteq\Gamma(\sh O_L)\otimes W.
    \]
    Indeed, we may identify 
    \[
      \Gamma(\sh O_L)=\GHom[_{\ger l_\ev}]0{\ger U(\ger l),\Gamma(\sh O_{L_0})},
    \]
    so we may consider superfunctions $f$ on $L$ as functions $f(u;g)$ of $u\in\ger U(\ger l)$ and $g\in L_0$. Consider the canonical extension $\mathbf a_u$ of the anti-homomorphism $x\longmapsto\mathbf a_x$ to $\ger U(\ger l)$ and $\mathbf a_g\defi \mathbf a\circ(j_g\times{\id})$ for $g\in L_0$. Then by \cite{ccf}*{Proposition 8.3.3}, we have for $w\in W$, $u\in\ger U(\ger l)$, and $g\in L_0$:
    \[
      \mathbf a^\sharp(w)(u;g)=\mathbf a_u\mathbf a_g^\sharp(w)=\mathsf S_\ger l(u)\cdot g^{-1}\cdot w\in W,
    \]
    with $\mathsf S_\ger l$ denoting the antipodal anti-automorphism of $\ger U(\ger l)$.

    Now, seeking a contradiction, assume that $p_\lambda^*=0$, so that 
    \[
      j_e^\sharp(\ell)=j_V^\sharp(\ell)(e)=0,\quad\forall \ell\in V_\lambda.
    \]
    Let $\ell\in L_V^\sharp$ be arbitrary. Then 
    \[
      \pi^\sharp j_V^\sharp(\ell)=\mathbf a_e^\sharp(\ell)=({\id}\otimes j_e)^\sharp\mathbf a^\sharp(\ell)\subseteq({\id}\otimes j_e)^\sharp(\Gamma(\sh O_L)\otimes V_\lambda)=0.
    \]
    Since $\pi$ is a surjective submersion, $\pi^\sharp$ is injective, so that $j_V^\sharp(\ell)=0$. But $j_V$ is an open embedding and $\ell$ is a superpolynomial, so $\ell=0$, contradiction! This proves the assertion. 
  \end{proof}
  
  \subsection{Determination of the spherical polynomials}
  \label{Sectstar}
  The Lie superalgebra $\ger g\cong\ger q(2n)$ carries a non-degenerate invariant odd supersymmetric bilinear form $b:\ger g\times\ger g\longrightarrow\cplxs$, the \Define{odd trace form}. In terms of the matrix realisation we have chosen in Equation \eqref{eq:gl-nonstd}, it is given by 
  \[
    b(x,x')\defi\frac12\tr(x\Phi x'),\quad x,x'\in\ger g
  \]
  or explicitly,
  \[
    \begin{split}
      b(x,x')&=\tr\Parens1{a\alpha'+\alpha a'+c\beta'+\beta c'+b\gamma'+\gamma b'+d\delta'+\delta d'},\\
      x&=
      \begin{Matrix}1
        a&\alpha&b&\beta\\\alpha&a&\beta&b\\
        \gamma&c&\delta&d\\c&\gamma&d&\delta
      \end{Matrix}\ ,\ x'=
      \begin{Matrix}1
        a'&\alpha'&b'&\beta'\\\alpha'&a'&\beta'&b'\\
        \gamma'&c'&\delta'&d\\c'&\gamma'&d'&\delta'
      \end{Matrix}.        
    \end{split}  
  \]
  The odd form restricts to an $\ger l$-invariant odd non-degenerate pairing of $V$ and $V^-$, allowing us to identify the $\ger l$-modules $V^*$ and $\Pi(V^-)$.

  We let $\ger t^*\defi\cplxs^n$. Then $\Pi(\ger t^*)$ may be identified with a subspace of $V^-$ via the odd map 
  \[
    \ger t^*\longrightarrow V^-:x\longmapsto 
    \begin{Matrix}1
      0&0&0&0\\
      0&0&0&0\\
      0&\diag(x)&0&0\\
      \diag(x)&0&0&0
    \end{Matrix}
  \]
  This determines an even injective linear map $\ger t^*\longrightarrow V^*$ and hence, a morphism $\aff(\ger t^*)\longrightarrow\aff(V^*)$. Thereby, we confer a meaning to $p|_{\ger t^*}$ for polynomials $p$ in $\mathscr P(V^*)=S(V)$, as the pullback along this morphism.
  
  Identifying $\Pi(V^-)$ with $V^*$, $p^*_\lambda|_{\ger t^*}$ is the ordinary $n$-variable polynomial given by 
  \begin{equation}\label{eq:nvarpol}
    p^*_\lambda|_{\ger t^*}(x_1,\dotsc,x_n)=p_\lambda^*
    \begin{Matrix}1
      0&0&0&0\\
      0&0&0&0\\
      0&\diag(x)&0&0\\
      \diag(x)&0&0&0
    \end{Matrix},
    \quad\forall x_j\in\cplxs.
  \end{equation}
  By \thmref{Prop}{polyact-ssmf}, $\smash{p_\lambda^*|_{\ger t^*}}$ is homogeneous of degree $\Abs0\lambda$ in $\mathscr P(V^*)$, so $p_\lambda^*|_{\ger t^*}$ is likewise homogeneous of degree $\Abs0\lambda$.
  
  \begin{Def}
    We call $p_\lambda^*\big|_{\ger t^*}$ the \Define{spherical polynomial} for $\lambda\in\Lambda^n_{>0}$.
  \end{Def}

  \medskip
  The remainder of this subsection is devoted to the proof of the following result.

  \begin{Th}[sphpol-q-hompart]
    Let $\lambda\in\Lambda^n_{>0}$ and $k\defi\Abs0\lambda$. Then the spherical polynomial $p^*_\lambda|_{\ger t^*}$ is the homogeneous part of degree $k$ of the eigenvalue polynomials $q_\lambda^*$, so 
    \[
      p^*_\lambda|_{\ger t^*}=\frac1{\lambda!}\prod_{i<j}\frac{\lambda_i-\lambda_j}{\lambda_i+\lambda_j}\,Q_\lambda.
    \]
  \end{Th}

  Let $\ger n^+$ be the subalgebra of $\ger l$ formed by the matrices
  \[
    \begin{Matrix}1
      A&0\\0&D
    \end{Matrix},\quad A=
    \begin{Matrix}1
      a&\alpha\\\alpha&a
    \end{Matrix},\quad D=
    \begin{Matrix}1
      d&\delta\\\delta&d
    \end{Matrix}
  \]
  where $a$ and $\alpha$ are strictly upper triangular and $d$ and $\delta$ are strictly lower triangular. Similarly, let $\ger a$ be the subspace of $\ger l$ formed by the matrices
  \[
    h_{a,\alpha}\defi\frac12
    \begin{Matrix}1
      A&0\\0&-A
    \end{Matrix},\quad A=
    \begin{Matrix}1
      a&\alpha\\\alpha&a
    \end{Matrix},
  \]
  where $a$ and $\alpha$ are diagonal. Then we have the vector space decomposition
  \[
    \ger l=\ger m\oplus\ger a\oplus\ger n^+.
  \]

  Although $\ger a_\ev$ is a subalgebra, $\ger a$ is not. We still consider the basis of $\ger a$ given by
  \[
    h_i\defi h_{e_i,0}=\frac12
    \begin{Matrix}[2]0
      E_{ii}&0&0&0\\0&E_{ii}&0&0\\0&0&-E_{ii}&0\\0&0&0&-E_{ii}
    \end{Matrix},\quad
    \eta_i\defi h_{0,e_i}=\frac12
    \begin{Matrix}[2]0
      0&E_{ii}&0&0\\E_{ii}&0&0&0\\0&0&0&-E_{ii}\\0&0&-E_{ii}&0
    \end{Matrix},
  \]
  where $i=1,\dotsc,n$.
  
  Let $\omega:S(\ger l)\longrightarrow\ger U(\ger l)$ denote the canonical supersymmetrisation map, given by 
  \[
    \omega(x_1\dotsm x_m)\defi\frac1{m!}\sum_{\sigma\in S_m}(-1)^{\#\Set0{(i,j)}{i<j,\sigma(i)>\sigma(j),\Abs0{x_i}=\Abs0{x_j}=\odd}}x_{\sigma(1)}\dotsm x_{\sigma(m)}
  \]
  \fa homogeneous $x_1,\dotsc, x_m\in\ger l$. 
  Then, by the Poincar\'e--Birkhoff--Witt theorem \cite{scheunert}*{Chapter I, \S~2.3, Corollary 1 to Theorem 1}, every $u\in\ger U(\ger l)$ can be written in a unique way as
  \begin{equation}\label{uu'usuba}
    u=u'+u_{\ger a},
  \end{equation}
  where
  \[
    u'\in \ger m\,\ger U(\ger l)+\ger U(\ger l)\ger n^+=\ger m\,\ger U(\ger l)\oplus\omega(S(\ger a\oplus\ger n^+))\ger n^+,
  \]
  and $u_{\ger a}$ can be written as a linear combination of monomials in the basis of $\ger a$, that is,
  \[
  u_{\ger a}=\sum_{J\in\mathcal J}c_Jh^J,
  \]
  where
  \[
    \mathcal J\defi\Set1{(j_1,\dotsc,j_n,j'_1,\dotsc,j'_n)}{\forall1\sle i\sle n:j_i\in\nats,j_i'=0,1},
  \]
  the constants $c_J\in\cplxs$, and 
  \[
    h^J\defi\prod_{i=1}^n h_i^{j_i}\prod_{i=1}^n \eta_i^{j'_i}.
  \]

  Recall that $z_\lambda$ is the unique element of $\mathcal Z(\ger q(n))$ such that $L(z_\lambda)=D_\lambda$. 

  \begin{Lem}[HCdecomposition]
    In the decomposition $z_\lambda=z_\lambda'+z_{\lambda, \ger a}$ of Equation \eqref{uu'usuba}, we have $z_\lambda\in\ger U(\ger a_\ev)$. In particular, there exists a unique polynomial $\gamma(z_\lambda)\in\cplxs[x_1,\dotsc,x_n]$ such that 
    \begin{equation}\label{eq:pbw-iwasawa}
      z_{\lambda,\ger a}=\gamma(z_\lambda)(h_1,\dotsc,h_n).
    \end{equation}
  \end{Lem} 

  \begin{proof}
    Let $\ger q(n)=\ger u^-\oplus\ger h\oplus\ger u^+$ be the standard triangular decomposition of $\ger q(n)$. Furthermore, let $\ger U(\ger q(n))^0$ denote the centralizer of $\ger h$ in $\ger U(\ger q(n))$. As usual, the Harish-Chandra projection gives rise to an $\mathrm{ad}_\ger h$-invariant direct sum decomposition
    \[
      \ger U(\ger q(n))^0=\Parens1{\ger u^-\,\ger U(\ger q(n))\cap \ger U(\ger q(n))\ger u^+}\oplus\ger U(\ger h).
    \] 
    Now write $z_\lambda\defi z_\lambda^0+z_{\lambda,\ger h}$, according to the latter decomposition. From $\mathrm{ad}_\ger h$-invariance of the decomposition, it follows that $[\ger h_\odd, z_{\lambda,\ger h}]=0$, and therefore that indeed $z_{\lambda,\ger h}\in\ger U(\ger h_\ev)$.  We can write $z_\lambda^0$ as a sum of monomials  of the form $x_1\dotsm x_r$, where $x_1,\dotsc, x_r\in\ger q(n)$,  such that $x_1\in\ger u^-$ and $x_r\in\ger u^+$. Next, by the natural embedding of $\ger U(\ger q(n))$ into $\ger U(\ger l)$ (as the left factor) we can consider the decomposition $z_\lambda\defi z_\lambda^0+z_{\lambda,\ger h}$ as one in $\ger U(\ger l)$, and clearly under this embedding the monomial $x_1\dotsm x_r$ is mapped to
    \[
      (x_1,0)\dotsm (x_r,0)\in\ger U(\ger l)\ger n^+.
    \]
    Similarly, since $z_{\lambda,\ger h}\in \ger U(\ger h_\ev)$, every monomial $x_1\dotsm x_r$ of $z_{\lambda,\ger h}$ can be written as 
    \[
      x_1\dotsm x_r=x_1'\dotsm x_r',
    \]
    where $x_i'\defi\frac{1}{2}\Parens1{(x_i,x_i)+(x_i,-x_i)}$, with 
    $(x_i,x_i)\in\ger m$ and $(x_i,-x_i)\in \ger a_\ev$. Since $\ger h_\ev$ is commutative, we obtain $z_{\lambda,\ger h}\equiv\frac{1}{2^r}\prod_{i=1}^r(x_i,-x_i)$ modulo $\ger m\,\ger U(\ger l)+\ger U(\ger l)\ger n^+$. The claim of the lemma now follows from  uniqueness of the decomposition $z_\lambda=z_{\smash\lambda}'+z_{\lambda,\ger a}$.
  \end{proof}

  \begin{Lem}[iwasawa-vanishing]
    Consider $e\in V$ as a linear polynomial in $\mathscr P(V^*)$. Then, for every non-negative integer $k\sge 0$, we have $\ger m\cdot e^k=0$ and 
    \[
      (\ger n^+\,\omega(S(\ger a\oplus\ger n^+))\cdot e^k)|_{\ger t^*}=0.
    \]
  \end{Lem}

  \begin{proof}
    We may assume that $k\sge1$. For $x\in\ger m$, we have 
    \[
      x\cdot e^k=k(x\cdot e)e^{k-1}=0,
    \]
    as $x\cdot e=0$. Next, Equations \eqref{eq:secondfactor-action} and \eqref{eq:firstfactor-action} imply that $\ger n^+$ acts by vector fields that lie in the span of 
    \[
      u_{pk}\frac\partial{\partial u_{p\ell}},\quad\xi_{pk}\frac\partial{\partial\xi_{p\ell}},\quad u_{pk}\frac\partial{\partial\xi_{p\ell}}, \quad\xi_{pk}\frac\partial{\partial u_{p\ell}},
    \]
    where $k>\ell$, whereas the action of $\ger a$ is given by vector fields of the same shape where now $k=\ell$. Let $\sh A^+$ be the subalgebra of $\sh P(V^*)$ generated by
    \[
      u_{\ell k},\quad\xi_{\ell k},\quad\forall k\sge\ell,
    \]
    and $\sh I^+$ the ideal of $\sh A^+$ generated by 
    \[
      u_{\ell k},\quad\xi_{\ell k},\quad\forall k>\ell.
    \]
    Then $\omega\Parens0{S(\ger a\oplus\ger n^+)}$ leaves $\sh A^+$ invariant and $\ger n^+$ maps $\sh A^+$ into $\sh I^+$. 

    On the other hand, $e$ is the sum of all the $u_{kk}$ and $\xi_{\ell\ell}$, and hence $e^k$ is contained in $\sh A^+$. Thus, $\ger n^+\,\omega(S(\ger a\oplus\ger n^+))\cdot e^k$ is contained in $\sh I^+$. As the restriction of $\sh I^+$ to $\ger t^*$ vanishes, this proves the claim. 
  \end{proof}

  For $\lambda\in\Lambda^n_{>0}$, let $D_\lambda^*\in\mathscr{PD}(V^*)$ be the operator dual to $D_\lambda$, \ie 
  \[
    \Dual0{D_\lambda^*(p)}q\defi\Dual0p{D_\lambda(q)},\quad\forall p\in\mathscr P(V^*),q\in\mathscr P(V).
  \]

  \begin{Prop}[sph-doexpr]
    \Fa $\lambda\in\Lambda^n_{>0}$ and $k\defi\Abs0\lambda$, we have 
    \begin{equation}\label{eq:sph-doexpr}
      p^*_\lambda=(-1)^kD_\lambda^*(e^k).
    \end{equation}
  \end{Prop}

  \begin{proof}
    The polynomial $e^k\in\mathscr P(V^*)$ is homogeneous of degree $k$. Thus, by \thmref{Prop}{polyact-ssmf} and \thmref{Cor}{ev-props}, we see that $D_\lambda^*(e^k)\in (V_\lambda)^*$. Furthermore, for every $\ell\in V_\lambda\subseteq\mathscr P(V)$, we have 
    \[
      \Dual1{D_\lambda^*(e^k)}\ell=c_\lambda(\lambda)\Dual0{e^k}{\ell}=k!\,j_V^\sharp(\ell)(e)=(-1)^k\Dual0{p_\lambda^*}\ell,
    \]
    in view of Equation \eqref{eq:sphpol-def}, and because $\ell$ is homogeneous of degree $k$. 
  \end{proof}

  \begin{proof}[\prfof{Th}{sphpol-q-hompart}]
    From the proof of \thmref{Prop}{ev-pol}, recall that there is an element $z_\lambda\in\sh Z(\ger q(n))\subseteq\sh Z(\ger l)$ (necessarily even) of order $\sle k$ \scth $L(z_\lambda)=D_\lambda$. Because of \thmref{Lem}{HCdecomposition}, we have 
    \[
      z_{\lambda,\ger a}=\gamma(z_\lambda)(h_1,\dotsc,h_n)\in\ger U(\ger a_\ev)=S(\ger a_\ev)
    \]
    for a unique polynomial $\gamma(z_\lambda)\in\cplxs[x_1,\dotsc,x_n]$. By Equation \eqref{eq:fpalg-action}, we have
    \[
      h_{a,0}\cdot e=a,\quad\forall a=\diag(a_1,\dotsc,a_n),
    \]
    so that 
    \begin{equation}\label{eq:hc-action}
      (z_{\lambda,\ger a}\cdot e^k)|_{\ger t^*}(x_1,\dotsc,x_n)=\gamma(z_\lambda)(x_1,\dotsc,x_n).
    \end{equation}
    Let $\mathsf S_{\ger l}:\ger U(\ger l)\longrightarrow \ger U(\ger l)$ denote the antipodal anti-automorphism of $\ger U(\ger l)$, defined by $\mathsf S_{\ger l}(x)\defi-x$ for $x\in\ger l$. Then
    \[
      \mathsf S_{\ger l}(z_{\lambda})\equiv (-1)^kz_{\lambda,\ger a}\pmod{\ger U(\ger l)\ger m\oplus{\ger n^+}\omega(S(\ger a\oplus\ger n^+))},
    \]
    as $\ger a_\ev$ is Abelian. Recall that we consider $\ger t^*=\cplxs^n$ as a subspace of $V^*$. Thus, for any $\mu\in\ger t^*$, we may form $\mu^k\in\mathscr P^k(V)=S^k(V^*)$. In these terms, we compute by the use of \thmref{Prop}{sph-doexpr} and \thmref{Lem}{iwasawa-vanishing}:
    \[
      \begin{split}
        p_\lambda^*|_{\ger t^*}(\mu)&=\frac1{k!}\Dual1{p_\lambda^*}{\mu^k}=\frac{(-1)^k}{k!}\Dual1{D_\lambda^*(e^k)}{\mu^k}\\
        &=\frac{(-1)^k}{k!}\Dual1{\mathsf S_{\ger l}(z_\lambda)(e^k)}{\mu^k}=\frac1{k!}\Dual1{z_{\lambda,\ger a}\cdot e^k}{\mu^k}.
      \end{split}
    \] 
    By Equation \eqref{eq:hc-action}, this is the value of the $k$-homogeneous part of $(z_{\lambda,\ger a}\cdot e^k)|_{\ger t^*}$ at the point $\mu\in\ger t^*$.

    Let $\mu\in\Lambda^n_{>0}$. By the highest weight theory for $\ger q(n)$ \cite{cwbook}*{\S~2.3}, there is a non-zero vector $v_\mu\in V_\mu$ \scth $\ger n^+\,v_\mu=0$ and 
    \[
      h_{a,0}\cdot v_\mu=\sum_{j=1}^n\mu_ja_j\cdot v_\mu,\quad\forall a=\diag(a_1,\dotsc,a_n).
    \]
    In view of Equation \eqref{eq:pbw-iwasawa}, and because $\eta_i^2=\tfrac12[\eta_i,\eta_i]\in\ger m$, we find that $V_\mu$ is generated as an $\ger m$-module by the vectors $v^i_\mu\defi\eta_iv_\mu$, $i=1,\dotsc,n$. Thus, there must be some $i$ \scth $\Dual0{p_\mu^*}{v^i_\mu}\neq0$. 

    As $\ger n^+\cdot v_\mu^i=0$ and $\ger m\cdot p_\mu^*=0$, we see that
    \[
      c_\lambda(\mu)\,\Dual0{p_\mu^*}{v^i_\mu}=\Dual0{p_\mu^*}{z_\lambda\cdot v^i_\mu}=\Dual0{p_\mu^*}{z_{\lambda,\ger a}\cdot v^i_\mu}=\gamma(z_\lambda)(\mu_1,\dotsc,\mu_n)\,\Dual0{p_\mu^*}{v^i_\mu}.
    \]
    We conclude
    \[
      q_\lambda^*(\mu)=c_\lambda(\mu)=\gamma(z_{\lambda,\ger a})(\mu_1,\dotsc,\mu_n)=(z_{\lambda,\ger a}\cdot e^k)|_{\ger t^*}(\mu_1,\dotsc,\mu_n),
    \]
    by Equation \eqref{eq:hc-action}, and with the convention from Equation \eqref{eq:nvarpol}. As $\Lambda^n_{>0}$ is Zariski dense in $\cplxs^n$, $p_\lambda^*|_{\ger t^*}$ is the $k$-homogeneous part of $q_\lambda^*$. The assertion now follows directly from \thmref{Th}{evpol-q}.
  \end{proof}

  In the next theorem, we use the notation introduced in Section \ref{SEcIntro}. In particular, recall that $\eta:\mathcal Z(\ger q(n))\longrightarrow \mathscr P(\ger h_\ev)$ denotes the Harish-Chandra homomorphism of Equation \eqref{HCeta}, and we identify the image of $\eta$ with the algebra of $Q$-symmetric polynomials in $n$ variables. Let $\mathsf S:\mathfrak U(\ger q(n))\longrightarrow \mathfrak U(\ger q(n))$ denote the canonical anti-automorphism of the enveloping algebra $\mathfrak U(\ger q(n))$ obtained by extending the map $\ger q(n)\longrightarrow\ger q(n):x\longmapsto -x$.

  \begin{Th}[cor1.3]
    For every $\lambda\in\Lambda_{>0}$, there exists a unique element $z_\lambda\in \mathcal Z(\ger q(n))$ such that $L(z_\lambda)=D_\lambda$.
    Furthermore,  
    \begin{equation}\label{*ezla}
      z_\lambda= r_\lambda\mathsf S(C_\lambda),
    \text{ where }
    r_\lambda\defi\frac{(-1)^{\lambda|}}{2^{\ell(\lambda)}(\lambda_1!\cdots\lambda_{\ell(\lambda)}!)^2}
    \prod_{1\sle i<j\sle \ell(\lambda)}
    \left(\frac{\lambda_i-\lambda_j}{\lambda_i+\lambda_j}\right)^2.
    \end{equation}
  \end{Th}  
    
  \begin{proof}
    Existence and uniqueness of $z_\lambda$ follow as in the proof of \thmref{Prop}{ev-pol}.
    Note that Nazarov considers a slightly different action of $\ger l$ on a polynomial space which decomposes as a direct sum of $\ger l$-modules of the form $2^{-\delta(\lambda)}F_\lambda\otimes F_\lambda$. Therefore, it follows from \cite{nazarov}*{Proposition 4.3} that $(-1)^{\Abs0\lambda}L(\mathsf S(C_\lambda))$ acts on $V_\lambda$ with the same scalar that $C_\lambda$ acts on Nazarov's $W_\lambda$ (see \cite{nazarov}*{Section 4}).
     
    For the proof of Equation \eqref{*ezla}, note first that the ``$\vrho$-shift'' for Sergeev's Harish-Chandra homomorphism is zero, and that $L(z_\lambda)$ acts on the modules $V_\mu\subseteq \mathscr P(V)$ by a scalar given by a $Q$-symmetric polynomial in $\mu_1,\dotsc,\mu_n$ of degree  at most $\Abs0\mu$. From \thmref{Prop}{PrptiQlQlstar} and \thmref{Cor}{ev-props}, it follows that $\eta(z_\lambda)=q_\lambda^*$. In particular,  the action of $L(z_\lambda)$ on $V_\mu$ is by given a polynomial in $\mu$ whose leading coefficient is $\smash{(Q_\lambda^*(\lambda))^{-1}}$. In addition, it follows from \cite{nazarov}*{Proposition 4.8} that the action of $(-1)^{\Abs0\lambda}L(\mathsf S(C_\lambda))$ is by a polynomial in $\mu$ with leading coefficient
    \[
      \lambda_1!\cdots\lambda_{\ell(\lambda)}!\prod_{1\sle i<j\sle \ell(\lambda)}\frac{\lambda_i+\lambda_j}{\lambda_i-\lambda_j}.
    \]    
     Equation \eqref{*ezla} now follows from a comparison of leading coefficients and the  formula for $Q_\lambda^*(\lambda)$ given in Ref.~\cite{ivanovv}.
  \end{proof}

\end{document}